\documentclass[10pt]{article}
\usepackage{amsmath}
\newtheorem{theorem}{Theorem}

\newtheorem{lemma}{Lemma}
\newtheorem{proposition}{Proposition}
\title{Function Theory for Laplace and Dirac-Hodge Operators in Hyperbolic Space}
\author{Yuying Qiao \thanks{Department of Mathematics, Hebei Normal University,
Shijiazhuang, P. R. China, Research supported by the National Science Foundation of China
(Mathematics Tianyuan Foundation, No A324610) and Hebei Province (105129),
E. Mail address {\tt yuyingqiao@163.com}} \and Swanhild Bernstein
\thanks{Institute of Mathematics and Physics, Bauhaus University Weimar, D-99421 Weimar, Germany,
E. Mail address {\tt swanhild.bernstein@fossi.uni-weimar.de}} \and Sirkka-Liisa Eriksson
\thanks{Department of Mathematics, University of Tampere, Tampere, Finland, Research supported by Academy of Finland, E. Mail address  {\tt Sirkka-Liisa.Eriksson@tut.fi}} \and John Ryan \thanks{Department of Mathematics, University of Arkansas, Fayetteville, AR 72701, USA, E.Mail address {\tt jryan@uark.edu}}}
\date{~}
\begin{document}
\maketitle
\begin{abstract}\ We develop basic properties of solutions to the Dirac-Hodge and Laplace equations in upper half space endowed with the hyperbolic metric. Solutions to the Dirac-Hodge equation are called hypermonogenic functions while solutions to this version of Laplace's equation are called hyperbolic harmonic functions. We introduce a Borel-Pompeiu formula for $C^{1}$ functions and a Green's formula for hyperbolic harmonic functions. Using a Cauchy Integral formula we are able to introduce Hardy spaces of solutions to the Dirac-Hodge equation. We also provide new arguments describing the conformal covariance of hypermonogenic functions and invariance of hyperbolic harmonic functions. We introduce intertwining operators for the Dirac-Hodge operator and hyperbolic Laplacian.
\end{abstract}


\section{Introduction}

\ Function theory for Dirac operators on manifolds have been developed in \cite{ca,cn,m1}. For particular types of manifolds this function theory has been developed in detail in \cite{lr,kr,kr1,kr2,r,vl} and elsewhere. In this paper we develop a detailed function theory associated to the Dirac-Hodge operator and Laplacian in upper half space, for $n>2$, endowed with the hyperbolic metric. Analysis of these operators have been developed over many years by many authors, particularly with respect to links to the Weinstein equation and its links to differential geometry and elasticity. See for instance \cite{a,al,cc,e1,e,el1,el2,el,el3,gs,h,hu,l1,l,w}.

\ Adapting the Cauchy Integral Formula for solutions to the Dirac-Hodge equation introduced in \cite{e} we also introduce a Green's formula for hyperbolic harmonic functions, Borel-Pompeiu formulas and other representation formulas. In particular we are able to study basic properties of Hardy spaces and Plemelj projection operators for hyper-surfaces in upper half space. We also introduce a Poisson integral formula to this setting. This thereby sets up the tools necessary for studying boundary value problems in this context. We also describe the conformal invariance of the operators introduced here and describe intertwining operators for these operators under actions of M\"{o}bius transformations preserving upper half space.

\section{Preliminaries}

\ Here we will consider upper half space $R^{n+}$ endowed with the hyperbolic metric $ds^{2}=\frac{dx_{1}^{2}+\ldots
+dx_{n}^{2}}{x_{n}^{2}}$. With respect to this metric one may consider the adjoint $\delta$ to the de Rham exterior
derivative $d$. Namely $\delta=\star d\star$, where $\star$ is the Hodge star map acting on sections in the alternating
bundle over $R^{n,+}$. The \emph{Dirac-Hodge operator} is the differential operator $d+\delta$ acting on differentiable
sections on the alternating algebra $\Lambda(R^{n,+})$. The square of $d+\delta$ is the Laplacian $d\delta+\delta d$
with respect to the hyperbolic or Poincar\/{e} metric. To better understand the Dirac-Hodge operator let us first
follow  \cite{l} and note that as an vector space the alternating or exterior algebra $\Lambda(R^{n})$ is isomorphic to
the Clifford algebra $Cl_{n}$ generated from $R^{n}$ with \emph{negative} definite inner product.  Namely let us
consider $R^{n}$ with orthogonal basis $e_{1},\ldots,e_{n}$. Then $Cl_{n}$ has as its basis

\[1, e_{1},\ldots,e_{n}, e_{1}e_{2},\ldots,e_{n-1}e_{n},\ldots,e_{1}\ldots,e_{n}\]
and $$e_{1}e_{j}+e_{j}e_{i}=-2\delta_{ij}.$$

\ Hence an arbitrary element of the basis may be written as $e_{A}=e_{\alpha_{1}}\ldots e_{\alpha_{h}}$,
here $A=\{\alpha_{1},\ldots,\alpha_{h} \}\subseteq\{1,\ldots,n\}$ and $1\le\alpha_{1}<\alpha_{2}<\ldots<\alpha_{h}\le n$.

\ We may express the Clifford algebra as $$Cl_{n}=Cl_{n-1}+Cl_{n-1}e_{n},$$ where $Cl_{n-1}$ is the Clifford algebra
generated from $R^{n-1}$ with orthonormal basis $e_{1},\ldots, e_{n-1}$. So if $A\in Cl_{n}$ there are unique elements
$B$ and $C\in Cl_{n-1}$ such that $A=B+Ce_{n}$. This gives rise to a pair of projection maps

\[P:Cl_{n}\rightarrow Cl_{n-1}:P(A)=B\]and \[Q:Cl_{n}\rightarrow Cl_{n-1}:Q(A)=C.\]

We will denote $-e_{n}Q(A)e_{n},\in Cl_{n-1}$, by $Q'(A)$.

\ The Dirac-Hodge operator, $d+\delta$ now retranslates in Clifford algebra notation as $D+\frac{n-2}{x_{n}}Q',$ where
$D=\Sigma_{j=1}^{n}e_{j}\frac{\partial}{\partial x_{j}}$ is the euclidean Dirac operator. So the Dirac-Hodge equation
is $$Df+\frac{n-2}{x_{n}}Q'(f)=0$$ where $f:U\rightarrow Cl_{n}$ is a differentiable function and $U$ is a domain in
$R^{n+}=\{x=x_{1}e_{1}+\ldots+x_{n}e_{n}:x_{n}>0\}$. See \cite{l} for more details. We shall abbreviate the Dirac-Hodge
equation to $Mf=0$. It may readily be determined that:

\begin{proposition}Suppose that $U$ is a domain in upper half space then the space of all solutions to the
Dirac-Hodge equation $\{f(x):x\in U$ and $Mf=0\}$ is a right module with respect to the algebra $Cl_{n-1}$.
\end{proposition}

Note, \cite{el1}, that if $U$ is a domain in upper half space and $h:U\rightarrow Cl_{n}$ is a $C^{2}$ function then

\[-M^{2}h=\triangle P(h)-\frac{n-2}{x_{n}}\frac{\partial P(h)}{\partial x_{n}}+
\left(\triangle Q(h)-\frac{n-2}{x_{n}}\frac{\partial Q(h)}{\partial x_{n}}+\frac{n-2}{x_{n}^{2}}Q(h)\right)e_{n}\]

where $\triangle$ is the euclidean Laplacian.

\ In \cite{a} it is noted for any real valued function $u(x)$ defined on the domain $U$ then
$$\triangle u-\frac{n-2}{x_{n}}\frac{\partial u}{\partial x_{n}}$$ is the Laplace formula
for upper half space with respect to the hyperbolic metric. We will denote this Laplacian by $\triangle_{R^{n,+}}$.
We will call a $Cl_{n-1}$ valued solution to the equation
$$\triangle h-\frac{n-2}{x_{n}}\frac{\partial u}{\partial x_{n}}=0$$
a \emph{hyperbolic harmonic function.} It follows that if $f$ is hypermonogenic and $C^{2}$ then $P(f)$ is hyperbolic
harmonic. Furthermore we shall denote the operator $$\triangle -\frac{n-2}{x_{n}}\frac{\partial}{\partial
x_{n}}+\frac{n-2}{x_{n}^{2}}$$ by $\triangle'_{R^{n,+}}$. The equations $\triangle_{R^{n,+}}u=0$ and
$\triangle'_{R^{n,+}}u=0$ are both examples of the \emph{Weinstein equation.} See for instance \cite{al,l1,w} for
details.

\ Returning to the Clifford algebra, we will need the anti-automorphism

\[\sim :Cl_{n}\rightarrow Cl_{n}:\sim e_{j_{1}}\ldots e_{j_{r}}=e_{j_{r}}\ldots e_{j_{1}}.\]

One usually writes $\tilde{A}$ for $\sim A$. Also for $A$, $B\in Cl_{n}$ one has, \cite{p},
$\widetilde{AB}=\tilde{B}\tilde{A}$. So if $f:U\rightarrow Cl_{n}$ satisfies $Mf=0$ then $\tilde{f}$ satisfies $fM=0$
where $fM=\Sigma_{j=1}^{n}\frac{\partial f}{\partial x_{j}}e_{j}+\frac{n-2}{x_{n}}Q'(f)$.

\ Following \cite{a1,v} one may express any M\"{o}bius transformation, $\phi(x)$, over $R^{n}\cup\{\infty\}$ as
$(ax+b)(cx+d)^{-1}$ where $a$, $b$, $c$ and $d$ are products of vectors from $R^{n}$ and $\tilde{a}c$, $\tilde{c}d$,
$\tilde{d}b$ and $\tilde{b}a\in R^{n}$. Moreover $a$, $b$, $c$ and $d$ are all products of vectors from $R^{n}$ and we
may assume that $\tilde{a}d-\tilde{b}c=\pm 1$. This gives rise to a covering group, $V(n)$, of the group of M\"{o}bius
transformations over $R^{n}\cup\{\infty\}$. We will be interested in the subgroup $V(n-1)$ that acts on $R^{n-1}$. The
group $V(n)$ is often called the Vahlen group. \ Following \cite{cw} for any four vectors $w_{1}$, $w_{2}$, $w_{3}$ and
$w_{4}\in R^{n}$ we define their \emph{cross ratio,} $$[w_{1},w_{2},w_{3},w_{4}],\quad \mbox{to be}\quad
(w_{1}-w_{4})^{-1}(w_{1}-w_{3})(w_{2}-w_{3})^{-1}(w_{2}-w_{4}).$$

\ Taking $A=a_{0}+\ldots +a_{1\ldots n}e_{1}\ldots e_{n}\in Cl_{n}$ we define the norm of $A$ to be, as usual,
$\|A\|=(a_{0}^{2}+\ldots+a_{1\ldots n}^{2})^{\frac{1}{2}}$. Using the conjugation antiautomorphism $-:Cl_{n}\rightarrow
Cl_{n}:-(e_{j_{1}}\ldots e_{j_{r}})=(-1)^{r}e_{j_{r}}\ldots e_{j_{1}}$ it may be seen that $\|A\|^{2}$ is the real part
of $A\overline{A}$, where $\overline{A}$ denotes the conjugate of $A$. It may be seen that if
$A=\underline{x}_{1}\ldots\underline{x}_{k}$ and each $\underline{x}_{j}\in R^{n}$ for $1\leq l\leq k$ then
$\|A\|^{2}=\|\underline{x}_{1}\|^{2}\ldots\|\underline{x}_{k}\|^{2}$. Each M\"{o}bius transformation $\psi(x)$,
$=(ax+b)(cx+d)^{-1}$ can be expressed as $ac^{-1}\pm (cx\tilde{c}+d\tilde{c})^{-1}$ whenever $c\ne 0$ and $\psi(x)=
\alpha ax\tilde{a}+bd^{-1}$ for some $\alpha\in R$ whenever $c=0$. Consequently:

\begin{lemma}
For each M\"{o}bius transformation $\psi$

\[\|[w_{1},w_{2},w_{3},w_{4}]\|=\|[\psi(w_{1}),\psi(w_{2}),\psi(w_{3}),\psi(w_{4})]\|.\]
\end{lemma}

\ This invariance of the norm of the cross ratio is also noted in \cite{a}.

\ The \emph{Cayley transformation} of upper half space $R^{n,+}$, $=\{x=x_{1}e_{1}+\ldots+x_{n}e_{n}:x_{n}>0\}$ to the
unit ball is given by $$C(x)=(e_{n}x+1)(x+e_{n})^{-1}=e_{n}(x-e_{n})(x+e_{n})^{-1}.$$  This transformation maps $e_{n}$
to the origin. If we wanted to adapt this transformation to a Cayley type M\"{o}bius transformation that maps upper
half space to the unit ball and maps a point $y$ in upper half space to the origin then one has the M\"{o}bius
transformation $$C(x,y)=e_{n}(x-y)(x-\hat{y})^{-1}$$ where $\hat{y}=y_{1}e_{1}+\ldots+y_{n-1}e_{n-1}-y_{n}e_{n}$. So
$\hat{y}$ is the reflection of $y$ about $R^{n-1}=span<e_{1},\ldots,e_{n-1}>$. Note that
$$\|C(x,y)\|=\frac{\|x-y\|}{\|x-\hat{y}\|}=\|[x,\hat{x},y,\hat{y}]\|^{\frac{1}{2}}.$$  Consequently we have the
following simple but important result.

\begin{lemma}
Suppose that $\psi\in V(n-1)$. Then $\psi(\hat{y})=\hat{\psi}(y)$ and
\[\|C(x,y)\|=\|C(\psi(x),\psi(y))\|.\]
\end{lemma}

\ As a consequence of this lemma one also has:

\begin{proposition}
Suppose $f:[0, \infty)\rightarrow Cl_{n-1}$ is an $L^{1}$ function and $\psi\in V(n-1)$ then
$$F(x,y)=\int_{0}^{\frac{\|x-y\|}{\|x-\hat{y}\|}}f(r)dr$$ is a well defined function on $R^{n,+}\times R^{n,+}$ and
$F(\psi(x),\psi(y))=F(x,y)$.
\end{proposition}

\section{Some Cauchy and Green's Integral Formulas}

\ Following \cite{a} let us first note that the hyperbolic Laplace equation on the unit ball in $R^{n}$ is

\[\triangle_{B(0,1)}u=\triangle u+\frac{2(n-2)r}{1-r^{2}}\frac{\partial u}{\partial r}=0.\]

\ Again following \cite{a} suppose now that $u(x)$ is a hyperbolic harmonic function depending only on $r=|x|$. First one obtains

\[\frac{\partial u}{\partial x_{i}}=u'(r)\frac{x_{i}}{r}\]

and

\[\frac{\partial ^2u}{\partial x_{i}\partial x_{j}}=u''(r)\frac{x_{i}x_{j}}{r^{2}}+
u'(r)\left(\frac{\delta_{ij}}{r}-\frac{x_{i}x_{j}}{r^{3}}\right).\]

Thus

\[\triangle u=u''+(n-1)\frac{u'}{r}.\]

\ As $u(x)$ is a hyperbolic harmonic function $u(r)$ will satisfy

\[u''+(n-1)\frac{u'}{r}+\frac{2(n-2)}{1-r^{2}}ru'=0.\]
If $u'\ne 0$ this can be written as
\[\frac{u''}{u'}+\frac{n-1}{r}+\frac{2(n-2)r}{1-r^{2}}=0.\]
or
\[\frac{d}{dr}[\log u'+(n-1)\log r-(n-2)\log(1-r^2)]=0\]
from which we conclude that
\[u'(r)\frac{r^{n-1}}{(1-r^{2})^{n-2}}=const.\]

\ This leads to the general solution

\[u(r)=a\int_{1}^{r}\frac{(1-t^{2})^{n-2}}{t^{n-1}}dt+b.\]

\ We see at once that no solution can stay finite for $r=0$. As a normalized solution we introduce
$$g(r)=\int\limits^1_r\frac{(1-t^2)^{n-2}}{t^{n-1}}dt.$$
 From Proposition 1 it now follows that the real valued function
 $$G(x,y)=\int_{\frac{\|x-y\|}{\|x-\hat{y}\|}}^{1}\frac{(1-t^{2})^{n-2}}{t^{n-1}}dt$$
  is a hyperbolic harmonic function. As $G(x,y)$ is real valued then trivially $Q(G(x,y))=0$. Consequently $MG(x,y)=DG(x,y)$. Therefore, \cite{l}, the function $p(x,y)=DG(x,y)$ is a vector valued hypermonogenic function. Here $M$ and $D$ are acting with respect to the $x$ variable. Following \cite{e, el} it may be noted that

\[DG(x,y)=\frac{(1-s^{2})^{n-2}}{s^{n-1}}\Big|_{\frac{\|x-y\|}{\|x-\hat{y}\|}}^{1}D\frac{\|x-y\|}{\|x-\hat{y}\|}\]
\[=\frac{(4x_{n}y_{n})^{n-2}}{\|x-y\|^{n-1}\|x-\hat{y}\|^{n-3}}\Sigma_{j=1}^{n}e_{j}\frac{\partial}{\partial x_{j}}\frac{\|x-y\|}{\|x-\hat{y}\|}\]
\[=\frac{(4x_{n}y_{n})^{n-2}}{\|x-y\|^{n-1}\|x-\hat{y}\|^{n-3}}
\left(\frac{x-y}{\|x-y\|\|x-\hat{y}\|}-\frac{(x-\hat{y})\|x-y\|}{\|x-\hat{y}\|^{3}}\right)\]
\[=(4x_{n}y_{n})^{n-2}\left(\frac{(x-y)^{-1}}{\|x-y\|^{n-2}\|x-\hat{y}\|^{n-2}}-
\frac{(x-\hat{y})^{-1}}{\|x-y\|^{n-2}\|x-\hat{y}\|^{n-2}}\right)\]
\[=x_{n}^{n-2}y_{n}^{n-1}\left(\frac{(x-y)}{\|x-y\|^{n}}e_{n}\frac{(x-\hat{y})}{\|x-\hat{y}\|^{n}}\right).\]

\ Suppose now that $U$ is a domain in upper half space and for two $C^{1}$ functions $f$ and $g$ defined on $U$ and
taking values in $Cl_{n}$ we consider the integral $\int_{S}g(x)\frac{n(x)}{x_{n}^{n-2}}f(x)d\sigma(x)$, where $S$
is a compact smooth hypersurface lying in $U$, $n(x)$ is the outer unit normal vector to $S$ at $x$ and $\sigma$ is the
Lebesgue surface measure of $S$. On assuming that $S$ is the boundary of a bounded subdomain $V$ of $U$ then
on applying Stokes' Theorem we obtain
\begin{multline*}
\int_{S}g(x)\frac{n(x)}{x_{n}^{n-2}}f(x)d\sigma(x)\\
=\int_{V}((g(x)D)\frac{1}{x_{n}^{n-2}}f(x)+g(x)\frac{1}{x_{n}^{n-2}}Df(x)-g(x)\frac{(n-2)}{x_{n}^{n-1}}e_{n}f(x)dx^{n}.
\end{multline*}
It follows that:

\begin{lemma} \cite{el} Suppose, $f$, $g$, $U$, $S$ and $V$ are as in the previous paragraph. Then
\[P\left(\int_{S} g(x)\frac{n(x)}{x_{n}^{n-2}}f(x)d\sigma(x)\right)=P\left(\int_{V}(g(x)M))f(x)+g(x)(Mf(x))\,
\frac{dx^{n}}{x_{n}^{n-2}}\right).\]
\end{lemma}

\ Consequently if $Mf=0$ and $gM=0$ we have the  version of Cauchy's Theorem established in \cite{el}. Namely $\int_{S}g(x)\frac{n(x)}{x_{n}^{n-2}}f(x)d\sigma(x)=0$. It may now be determined that for each $y\in V$
\[P(f(y))=\frac{2^{n-2}}{\omega_{n}}P\left(\int_{S}p(x,y)\frac{n(x)}{x_{n}^{n-2}}f(x)d\sigma(x)\right).\]
This is the Cauchy integral formula arising in \cite{el}. It is an easy consequence of this integral formula and the previous lemma to obtain:

\begin{theorem}{\bf{(Borel Pompeiu Theorem)}} Suppose that $f:U\rightarrow Cl_{n}$ is a $C^{1}$ function and that $U$ is a bounded open subset of upper half space with $C^{1}$ compact boundary lying in upper half space. Suppose also that $f$ has a continuous extension to the boundary of $U$. Then for each $y\in U$

\[P(f(y))=\frac{1}{\omega_{n}}P\left(\int_{\partial U}p(x,y)\frac{n(x)}{x_{n}^{n-2}}f(x)d\sigma(x)+\int_{U}p(x,y)
(Mf(x))\frac{dx^{n}}{x_{n}^{n-2}}\right).\]\end{theorem}

\ Clearly if $f(y)\in Cl_{n-1}$ then this integral would give $f(y)$.

\ It follows from this integral formula that if $\phi$ is a $C^{\infty}$ function with values in $Cl_{n-1}$ and with compact support on upper half space then

\[\phi(y)=\frac{1}{\omega_{n}}\int_{R^{n,+}}p(x,y)(M\phi(x))\frac{dx^{n}}{x_{n}^{n-2}}\]
for each $y\in R^{n,+}$.\ We also have as a consequence of Lemma 3 the following version of Green's Representation Formula for hyperbolic harmonic functions.

\begin{theorem}{\bf{(Green's Formula)}}Suppose that $U$ is a domain in upper half space and that $h:U\rightarrow Cl_{n-1}$ is a hyperbolic harmonic function. Then for $S$ a piecewise $C^{1}$, compact surface lying in $U$ and bounding a bounded subdomain $V$ of $U$

\[h(y)=\frac{1}{\omega_{n}}P\left(\int_{S}G(x,y)\frac{n(x)}{x_{n}^{n-2}}\Big(Mh(x))
-p(x,y)\frac{n(x)}{x_{n}^{n-2}}h(x)\Big)d\sigma(x)\right)\]

for each $y\in V$.
\end{theorem}

\ Stokes' Theorem also gives that if $\phi$ is $Cl_{n-1}$ valued, $C^{\infty}$, is defined on upper half space and has compact support then for each $y\in R^{n,+}$

\begin{equation}
\phi(y)=\frac{1}{\omega_{n}}\int_{R^{n,+}}G(x,y)\big(\triangle_{R^{n,+}}\phi(x)\big)\frac{dx^{n}}{x_{n}^{n-2}}.
\end{equation}

\ Now let us consider $D_{y}G(x,y)$ where $D_{y}=\Sigma_{j=1}^{n}e_{j}\frac{\partial}{\partial y_{j}}$. As $\|x-\hat{y}\|=\|y-\hat{x}\|$ then $G(x,y)$ is hyperbolic harmonic in both the variables $x$ and $y$, and

\[D_{y}G(x,y)=D_{y}\int_{\frac{\|y-x\|}{\|y-\hat{x}\|}}^{1}\frac{(1-s^{2})^{n-2}}{s^{n-1}}ds\]
\[=4x_{n}^{n-2}y_{n}^{n-2}\left(\frac{(y-x)^{-1}}{\|x-y\|^{n-2}\|y-\hat{x}\|^{n-2}}
-\frac{(y-\hat{x})^{-1}}{\|x-y\|^{n-2}\|x-y\|^{n-2}}\right)=h(x,y)\] is hypermonogenic in the variable $y$.

\ Let $M_{y}$ denote the Dirac-Hodge operator with respect to the variable $y$ and let $\triangle_{R^{n,+},y}$ denote the hyperbolic Laplacian with respect to the variable $y$.

\begin{theorem}
Suppose that $\psi$ is a $Cl_{n-1}$ valued, $C^{\infty}$ function with compact support on upper half space. Then

\[P\left(M_{y}\left(\frac{1}{\omega_{n}}\int_{R^{n,+}}h(x,y)\psi(x)\frac{dx^{n}}{x_{n}^{n-2}}\right)\right)=\psi(y).\]
\end{theorem}

\ Now consider

\[\triangle_{R^{n,+},y}\left(\frac{1}{\omega_{n}}\int_{R^{n,+}}G(x,y)\psi(x)\frac{dx^{n}}{x_{n}^{n-2}}\right).\]
This is equal to
\[\frac{1}{\omega_{n}}P\left(M_{y}\left(D\int_{R^{n,+}}G(x,y)\psi(x)\frac{dx^{n}}{x_{n}^{n-2}}\right)\right),\]
which in turn is equal to
\[\frac{1}{\omega_{n}}P\left(M_{y}\left(\int_{R^{n,+}}h(x,y)\psi(x)\frac{dx_{n}}{x_{n}^{n-2}}\right)\right).\]
By Theorem 3 this evaluates to $\psi(y)$. So we have established:

\begin{theorem}
Suppose $\psi$ is as in Theorem 3 then

\[\triangle_{R^{n,+},y}\left(\frac{1}{\omega_{n}}\int_{R^{n,+}}G(x,y)\psi(x)\frac{dx^{n}}{x_{n}^{n-2}}\right)=\psi(y).\]
\end{theorem}

\ In \cite{e} the kernel $$q(x,y)=DH(x,y)=\frac{1}{2(n-2)}D\frac{1}{\|x-y\|^{n-2}\|x-\hat{y}\|^{n-2}},$$ where
$$H(x,y)=\frac{1}{(n-2)\|x-y\|^{n-2}\|x-\hat{y}\|^{n-2}}$$ is introduced. In \cite{e} it is shown that  the kernel
$q(x,y)$ is the Cauchy kernel for the $Q$ part of a Cauchy Integral Formula for hypermonogenic functions. So from
\cite{e} we have

\[f(y)=P(f(y))+Q(f(y))e_{n}=\]
\[\frac{2^{n-1}y_{n}^{n-2}}{\omega_{n}}\left(P\left(\int_{\partial U}r(x,y)\frac{n(x)}{x_{n}^{n-2}}f(x))d\sigma(x)\right)
-Q\left(\int_{\partial U}q(x,y)n(x)f(x)d\sigma(x)\right)e_{n}\right)\] where $r(x,y)=y_{n}^{-n+2}p(x,y)$.

\ Again as a consequence of Stokes' Theorem we have:

\begin{theorem}
Suppose that $\phi$ be a $Cl_{n}$ valued $C^{1}$ function defined on a bounded domain $U\subset R^{n,+}$,
with piecewise smooth boundary, and $\phi$ has a continuous extension to the closure of $U$. Then for each $y\in U$

\[Q(\phi(y))=\frac{2^{n-2}y_{n}^{n-2}}{\omega_{n}}Q\left(\int_{\partial U}q(x,y)n(x)\phi(x)d\sigma(x)
-\int_{U}q(x,y)(M\phi(x))dx^{n}\right).\]
\end{theorem}

\ It follows immediately that if $\phi$ has compact support then

\[Q(\phi(y))=\frac{2^{n-2}y_{n}^{n-2}}{\omega_{n}}Q\left(\int_{R^{n,+}}q(x,y)\big(M\phi(x)\big)dx^{n}\right).\]

\ Furthermore it may readily be determined that:

\begin{theorem}{\bf{(Green's Formula:)}}Suppose that $u:U\rightarrow Cl_{n-1}e_{n}$ is a solution of the equation $\triangle'_{R^{n,+}}u=0$, and $U$ is as in Theorem 5. Then for each $y\in U$ we have

\[u(y)=\frac{2^{n-2}y_{n}^{n-2}}{\omega_{n}}Q\left(\int_{\partial U} H(x,y)n(x)\big(Mu(x)\big)-q(x,y)n(x)u(x)
d\sigma(x)\right).\]
\end{theorem}

\ In particular if $u$ is a real valued function satisfying $\triangle'_{R^{n,+}}u=0$ then

\[u(y)=-e_{n}\frac{2^{n-2}y_{n}^{n-2}}{\omega_{n}}\int_{\partial U}H(x,y)n(x)M\big(e_{n}u(x)\big)-q(x,y)n(x)e_{n}u(x)d\sigma(x).\]

\ By similar arguments to those used before we also have:

\begin{theorem}
Suppose the $U$ is as in Theorem 5 and $u:U\rightarrow Cl_{n-1}e_{n}$ is a $C^{2}$ function then
\begin{eqnarray*}
u(y)=\frac{y_{n}^{n-2}}{\omega_{n}}Q\left(\int_{\partial U}H(x,y)n(x)\big(Mu(x)\big)-q(x,y)n(x)u(x)d\sigma(x) \right. \\
-\left.\int_{U}H(x,y)\big(\triangle'_{R^{n,+}}u(x)\big)dx^{n}\right).
\end{eqnarray*}
\end{theorem}

\ Consequently if $u$ has compact support then
\[u(y)=\frac{y_{n}^{n-2}}{\omega_{n}}\int_{R^{n,+}}H(x,y)\big(\triangle'_{R^{n,+}}u(x)\big)dx^{n}.\]

\ In \cite{l1} within Lemma 2.1 it is shown that if $\phi(x)$ is a solution to
$$\triangle\phi(x)-\frac{n-2}{x_{n}}\frac{\partial\phi(x)}{\partial x_{n}}+\frac{n-2}{x_{n}^{2}}\phi(x)=0$$ then
$\theta(x)=x_{n}^{n-2}\phi(x)$ is a solution to the equation
$$\triangle\theta(x)-\frac{n-2}{x_{n}}\frac{\partial\theta(x)}{\partial x_{n}}+\frac{n-2}{x_{n}}\theta(x)=0.$$ As
$\|\hat{y}-x\|=\|y-\hat{x}\|$ it follows that $y^{n-2}H(x,y)$ is hyperbolic harmonic in the $y$ variable. So by simple
adaptations of standard arguments we also have

\begin{proposition}
Suppose $u:R^{n,+}\rightarrow R$ is a $C^{2}$ function with compact support. Then

\[u(y)=\triangle'_{R^{n,+}}\left(\frac{1}{\omega_{n}}y_{n}^{n-2}\int_{R^{n,+}}H(x,y)u(x)dx^{n}\right)\]
for each $y\in R^{n,+}$.
\end{proposition}

\ Now for any $A\in Cl_{n}$, $P(A)=\frac{1}{2}(A+\hat{A})$ where $\hat {A}=B-Ce_{n}$ with $B$ and $C\in Cl_{n-1}$.
Moreover, $Q(A)=\frac{-1}{2}(A-\hat{A})e_{n}$ and for any elements $X$ and $Y\in Cl_{n}$ it is straightforward to
determine that $\widehat{XY}=\hat{X}\hat{Y}$. Using these observations it is noted in \cite{e} that the previous
integral becomes
\begin{eqnarray*}
\frac{1}{\omega_{n}}2^{n-1}y_{n}^{n-2}\left(\int_{\partial U}\frac{1}{2}\Big(r(x,y)n(x)\frac{n(x)}{x_{n}^{n-2}}f(x)
+\hat{r}(x,y)\frac{\hat{n}(x)}{x_{n}^{n-2}}\hat{f}(x)\Big)d\sigma(x) \right. \\
-\left.\int_{\partial U}\frac{e_{n}}{2}\Big(q(x,y)n(x)f(x)-\hat{q}(x,y)\hat{n}(x)\hat{f}(x)\Big)d\sigma(x)\right).
\end{eqnarray*}

\ In \cite{e} it is shown that this expression simplifies to
\begin{eqnarray*}
f(y)=\frac{2^{n-1}y_{n}^{n-2}}{\omega_{n}}\left(\int_{\partial
K}\frac{(x-y)^{-1}n(x)f(x)}{\|x-y\|^{n-2}\|x-\hat{y}\|^{n-2}}d\sigma(x) \right. \\
-\left. \int_{\partial
K}\frac{(\hat{x}-y)^{-1}\hat{n}(x)\hat{f}(x)}{\|x-y\|^{n-2}\|\hat{x}-y\|^{n-2}}d\sigma(x)\right).
\end{eqnarray*}

\ If we write $E(x,y)$ for $\frac{(x-y)^{-1}}{\|x-y\|^{n-2}\|x-\hat{y}\|^{n-2}}$ and $F(x,y)$ for $\frac{(\hat{x}-y)^{-1}}{\|x-y\|^{n-2}\|\hat{x}-y\|^{n-2}}$ then this integral formula simplifies to

\[f(y)=\frac{2^{n-1}y_{n}^{n-2}}{\omega_{n}}\int_{S}\Big(E(x,y)n(x)f(x)-F(x,y)\hat{n}(x)\hat{f}(x)\Big)d\sigma(x).\]

\section{Plemelj Projection Operators and Hardy Spaces of Hypermonogenic Functions}

\ First let us note that as $y_{n}$ tends to infinity then $y_{n}^{n-2}E(x,y)$ and $y_{n}^{n-2}F(x,y)$ both tend to zero for fixed $x$. Also as $y_{n}$ tends to zero then both $y_{n}^{n-2}E(x,y)$ and $y_{n}^{n-2}F(x,y)$ tend to zero for fixed $x$.

\begin{proposition}
Suppose that $C\in Cl_{n}$ is a constant and $S$ is a compact, $C^{2}$ surface lying in upper half space.
Suppose further that $S$ is the boundary of a bounded domain $U$ in $R^{n,+}$. If $y(t)$ is a $C^{1}$ path
in $U^{+}$ with nontangential limit $y(1)=y\in S$ then
\[\lim_{t\rightarrow 1}\frac{2^{n-2}y(t)_{n}^{n-2}}{\omega_{n}}\int_{S}\Big(E(x,y)n(x)C-
F(x,y)\hat{n}(x)\hat{C}\Big)d\sigma(x)=\]
\[\frac{1}{2}C+\frac{2^{n-2}y_{n}^{n-2}}{\omega_{n}}PV\int_{S}\Big(E(x,y)n(x)C-F(x,y)\hat{n}(x)\hat{C}\Big)d\sigma(x).\]
\end{proposition}
{\bf{Proof:}} Given that
$$\lim_{x\rightarrow y(1)}\frac{2^{n-2}y(1)_{n}^{n-2}}{\|\hat{x}-y(1)\|^{n-2}}=1$$
then as $S$ is compact it follows from the Mean Value Theorem that given $\epsilon>0$ then for all $x\in S$ such that
$\|x-y(1)\|<1$ we have
$$\left|\frac{2^{n-2}y(1)_{n}^{n-2}}{\|\hat{x}-y(1)\|^{n-2}}-1\right|<C'\|x-y(1)\|$$
and $C'\geq 0$. Let $S_{\epsilon}(y)=\{x\in S:\|x-y\|<\epsilon\}$. It follows that
\begin{align*}
\frac{2^{n-2}y(t)_{n}^{n-2}}{\omega_{n}}&\int_{S_{\epsilon}(y)}E(x-y)n(x)Cd\sigma(x)\\
=\frac{1}{\omega_{n}}&\int_{S}\frac{(x-y(t))}{\|x-y(t)\|^{n}}n(x)C\,d\sigma(x)\\
&+\frac{1}{\omega_{n}}\int_{S_{\epsilon}(y)}\left(\frac{2^{n-2}y(t)_{n}^{n-2}}{\|\hat{x}-y(t)\|^{n-2}}-1\right)
\frac{(x-y)}{\|x-y(t)\|^{n}}n(x)Cd\sigma(x).
\end{align*}
It follows from the usual calculations, see \cite{i}, for Plemelj formulas in Clifford analysis that
\[\lim_{t\rightarrow 1}\frac{1}{\omega_{n}}\int_{S_{\epsilon}(y)}\frac{\big(x-y(t)\big)}{\|x-y(t)\|^{n}}n(x)Cd\sigma(x)=\frac{1}{2}C\]
and the other term tends to zero at $t$ tends to $1$. The result now follows. Q.E.D.

\ Although in the previous proposition we assumed that the surface is $C^{2}$ one can also prove this result for surfaces that are strongly Lipschitz. These are surfaces that are locally Lipschitz graphs and whose Lipschitz constants are uniformly bounded.

\ One also readily has the following important technical result.

\begin{lemma}
For $x\in R^{n,+}$ and fixed $y\in R^{n,+}$ with $\|x-y\|>2y_{n}$ then $\|E(x,y)\|<\frac{C}{\|x-y\|^{2n-2}}$ and $\|F(x,y)\|<\frac{C}{\|x-y\|^{2n-2}}$ for some $C\in R^{+}$.
\end{lemma}

\ Using this lemma one can adapt arguments developed in \cite{glq,m} and elsewhere to deduce:

\begin{theorem}
Suppose that $\Sigma$ is a Lipschitz graph lying in upper half space and the minimal distance between $\Sigma$ and
the boundary, $R^{n-1}$, of upper half space is greater than zero then the singular integral operator $T_{\Sigma}$
defined by
\[\frac{2^{n-2}}{\omega_{n}}PV\int_{\Sigma}y_{n}^{n-2}\left(E(x,y)n(x)\phi(x)-F(x,y)\hat{n}(x)\hat{\phi}(x)\right)d\sigma(x)\]
is $L^{p}$ bounded for $1<p<\infty$.
\end{theorem}

\ Clearly this result also holds if we replace the Lipschitz graph $\Sigma$ by a compact, strongly Lipschitz surface $S$. In this case the operator $T_{\Sigma}$ is replaced by its analogue $T_{S}$.

\ This result enables us to establish the analogues of Plemelj formulas in the present context.

\begin{theorem}
Suppose that $S$ is a compact, strongly Lipschitz surface lying in upper half space. Suppose also that $S$ is the
boundary of a bounded domain $U^{+}$ and an exterior domain $U^{-}\subset R^{N,+}$. Then for each function
$\phi\in L^{p}(S)$ for $1<p<\infty$ or a path $y_{\pm}(t)\in U^{\pm}$ with nontangential limit $y(1)=y\in S$ we have
\[\lim_{t\rightarrow 1}\frac{2^{n-2}y_{\pm}(t)^{n-2}}{\omega_{n}}\int_{S}\left(E\big(x,y(t)\big)n(x)\phi(x)-F\big(x,y(t)\big)\hat{n}(x)\hat{\phi}(x)\right)d\sigma(x)\]
\[=\pm\frac{1}{2}\phi(y)+\frac{2^{n-2}}{\omega_{n}}PV\int_{S}y_{n}^{n-2}\left(E(x,y)n(x)\phi(x)-F(x,y)\hat{n}(x)\hat{\phi}(x)\right)d\sigma(x)\]
for almost all $y\in S$.
\end{theorem}

\ A minor adaptation of the proof of Theorem 17 in \cite{e} tells us the following:

\begin{theorem} Suppose $S$ is a Lipschitz surface lying in the closure of upper half space and $\phi\in L^{p}(S)$
for some $p\in(1,\infty)$ then the integral

\[\frac{2^{n-2}y_{n}^{n-2}}{\omega_{n}}\int_{S}\left(E(x,y)n(x)\phi(x)-F(x,y)\hat{n}(x)\hat{\phi}(x)\right)d\sigma(x)\]
defines a left hypermonogenic function $f(y)$ on $R^{n,+}\backslash S$.
\end{theorem}

\ As $\lim_{y_{n}\rightarrow\infty}y_{n}^{n-2}E(x,y)=0$ and $\lim_{y_{n}\rightarrow\infty}y_{n}^{n-2}F(x,y)=0$ for each $x\in S$ and
$\lim_{y_{n}\rightarrow 0}y_{n}^{n-2}E(x,y)=\lim_{y_{n}\rightarrow 0}y_{n}^{n-2}F(x,y)=0$ for each $x\in S$ then
$lim_{y_{n}\rightarrow\infty}f(y)=\lim_{y_{n}\rightarrow 0}f(y)=0$. It now follows that the operators
$$\frac{1}{2}I\pm T_{S}:L^{p}(S)\rightarrow L^{p}(S)$$ are projection operators with images the Hardy spaces
\begin{eqnarray*}
H^{p}(U^{\pm})=\{f:U^{\pm}\rightarrow Cl_{n}: f \ \mbox{is left hypermonogenic and}\\
\mbox{nontangentially approaches some element in}\ L^{p}(S)\}.
\end{eqnarray*}

\ Consequently

\[L^{p}(S)=H^{p}(U^{+})\oplus H^{p}(U^{-}).\]

\ The operators $\frac{1}{2}I\pm T_{S}$ are generalizations of the Plemelj projection operators to the context of
hypermonogenic functions. As in the euclidean case these operators are projection operators, or mutually annihilating
idempotents. Let us denote the operator $\frac{1}{2}I+T_{S}$ by $\cal{H}_{S}$. We may introduce the \emph{Kerzman-Stein
operator} $A_{S}=\cal{H}_{S}-\cal{H}_{S}^{\star}$, where $\cal{H}_{S}^{\star}$ is the adjoint of $\cal{H}_{S}$. In
particular if $\phi\in L^{2}(S)$ then
\begin{multline*}
A_{S}(\phi)=\frac{2^{n-2}y_{n}^{n-2}}{\omega_{n}}\left(\int_{\partial K}
\big(E(x,y)n(x)-n(x)E(x,y)\big)\phi(x)\right. \\
-\left. \big(F(x,y)\hat{n}(x)+\hat{n}(x)F(x,y)\big)\hat{\phi}(x)d\sigma(x)\right) .
\end{multline*}
\ Let us now turn to consider the case where $S=R^{n,+}$. We begin with:

\begin{theorem}
Suppose $\phi\in L^{p}(R^{n-1})$ for some $p\in(1,\infty)$ then for $y(t)=y'+y_{t}e_{n}$, where $y'\in R^{n-1}$
and $y_{n}(t)>0$,

\[\lim_{t\rightarrow 0}\frac{2^{n-2}y_{n}(t)^{n-2}}{\omega_{n}}\left(\int_{R^{n-1}}E(x,y(t))e_{n}\phi(x)+F(x,y(t))e_{n}
\hat{\phi}(x)dx^{n-1}\right)=P(\phi(y'))\] almost everywhere.
\end{theorem}
{\bf{Proof:}} Without loss of generality we may assume that $y'=\underline{0}$. Let us assume that $\phi$ is
$C^{\infty}$ and has compact support. Now for any $\epsilon>0$

\[\lim_{t \rightarrow 0}\frac{2^{n-2}y_{n}(t)^{n-2}}{\omega_{n}}\left(\int_{R^{n-1}\backslash B(0,\epsilon)}E(x,y_{n}e_{n}
(t))e_{n}\phi(x)+F(x,y_{n}(t)e_{n})e_{n}\hat{\phi}(x)dx^{n-1}\right)\] is equal to zero.

\ Further by arguments similar to those used to establish Proposition 4 we have:

\[\lim_{\epsilon\rightarrow 0, t\rightarrow 0}\frac{2^{n-2}y_{n}(t)^{n-2}}{\omega_{n}}\int_{B(0,\epsilon)}
E(x,y_{n}(t)e_{n})e_{n}\phi(x)dx^{n-1}=\frac{1}{2}\phi(0).\] As $\hat{x}=x$ for each $x\in R^{n-1}$ then similarly

\[\lim_{\epsilon\rightarrow 0,t\rightarrow 0}\frac{2^{n-2}y_{n}(t)^{n-2}}{\omega_{n}}\int_{B(0,\epsilon)}
F(x,y_{n}(t)e_{n})e_{n}\hat{\phi}(x)dx^{n-1}=\frac{1}{2}\hat{\phi}(y').\] A standard density argument now reveals the
result for all $\phi\in L^{p}(R^{n-1})$. Q.E.D.

\ Theorem 11 tells us that in the special case where $S=R^{n-1}$ we may solve a Dirichlet problem for $Cl_{n-1}$ valued $L^{p}$ boundary data and for the Dirac-Hodge equation as opposed to the hyperbolic Laplace equation. This is in contrast to the euclidean analogues where one obtains a Plemelj formula for such data..

\ Suppose now that $\phi(x)$ is a real valued, $L^{p}$ function defined on $R^{n-1}$, with $1<p<\infty$. Then on restricting to the real part of our previous integral we have the following Poisson integral

\[F(y)=\frac{2^{n-2}y_{n}^{n-1}}{\omega_{n}}\left(\int_{R^{n}}\frac{\phi(x)}{\|x-y\|^{n}\|x-\hat{y}\|^{n-2}}+
\frac{\phi(x)}{\|x-y\|^{n-2}\|x-\hat{y}\|^{n}}d\sigma(x)\right),\] which defines a hyperbolic harmonic function on
upper half space with boundary value $\phi$.

\ The material developed in this section enables one to tackle boundary values problems for the hyperbolic harmonic
equation and for the equation $\triangle'_{R^{n,+}}u=0$. This includes problems like the Dirichlet and Neumann problems.
One may adapt arguments given in \cite{Mc,m} to the context described here and solve such boundary value problems for
hyperharmonic functions. This will be done elsewhere.

\section{Representation Theorems}

\ We begin with:

\begin{theorem}{\bf{(Borel-Pompeiu Formula)}}Let $K\subset R^{n,+}$ be a bounded region with smooth boundary in $R^{n,+}$, Suppose also that $f:K\rightarrow Cl_{n}$ is a $C^{1}$ function on $K$ with a continuous extension to the closure of $K$. Then for $y\in K$ we have

\[f(y)=\frac{(2y_{n})^{n-1}}{\omega_{n}}\int_{\partial K}P\Big(p(x,y)\frac{n(x)}{x_{n}^{n-2}}f(x)\Big)d\sigma(x)
+\frac{1}{y_{n}}Q\big(q(x,y)n(x)f(x)\big)d\sigma(x)e_{n}\]
\[-\frac{(2y_{n})^{n-1}}{\omega_{n}}\int_{K}\left[P\big(p(x,y)Mf\big)\frac{1}{x_{n}^{n-1}}+Q\big(q(x,y)Mf\big)
\frac{e_{n}}{y_{n}}\right]dx^{n}\] or
\[f(y)=\frac{(2y_{n})^{n-2}}{\omega_{n}} \int_{\partial K}E(x,y)n(x)f(x)d\sigma(x)\]
\[-\int_{\partial K}F(x,y)\hat{n}(x)\hat{f}(x)d\sigma(x)
-\frac{(2y_{n})^{n-2}}{\omega_{n}}\int_{K}\big(E(x,y)Mf-F(x,y)\widehat{Mf}\big)dx^{n}.\]
\end{theorem}
{\bf{Proof:}} Consider a sphere $U(y,\delta)\subset K$ with center at $y$ and radius $\delta>0$ then we have

\[\int_{\partial K}P\left(p(x,y)\frac{n(x)}{x_{n}^{n-2}}f(x)d\sigma(x)\right)+
\frac{1}{y_{n}}\int_{\partial K}Q\left(q(x,y)n(x)f(x)d\sigma(x)\right)e_{n} \]
\[=\int_{\partial K\backslash\partial U(y,\delta)}P\left(p(x,y)\frac{n(x)}{x_{n}^{n-2}}f(x)d\sigma(x)\right)\]
\[+\frac{1}{y_{n}}\int_{\partial K\backslash\partial U(y,\delta)}Q\left(q(x,y)n(x)f(x)d\sigma(x)\right)e_{n} \]
\[+\int_{\partial U(y,\delta)}P\left(p(x,y)\frac{n(x)}{x_{n}^{n-2}}f(x)d\sigma(x)\right)
+\frac{1}{y_{n}}\int_{\partial U(y,\delta)}Q\left(q(x,y)n(x)f(x)d\sigma(x)\right)e_{n} \]
\[=\int_{K\backslash U(y,\delta)}P\left(p(x,y)Mf\right)\frac{1}{x_{n}^{n-1}}+\frac{1}{y_{n}}Q\left(q(x,y)Mf\right)dx^{n}\]
\[+ \int_{\partial U(y,\delta)}P\left(p(x,y)\frac{n(x)}{x_{n}^{n-2}}f(x)d\sigma(x)\right)
+\frac{1}{y_{n}}\int_{\partial U(y,\delta)}Q\left(q(x,y)n(x)f(x)d\sigma(x)\right)e_{n}\]

\ When $\delta$ tends to $0$ then

\[\int_{\partial U(y,\delta)}P\left(p(x,y)\frac{n(x)}{x_{n}^{n-2}}f(x)d\sigma(x)\right)+\frac{1}{y_{n}}\int_{\partial U(y,\delta)}Q(q(x,y)n(x)f(x)d\sigma(x))e_{n} \]
tends to $\frac{\omega_{n}}{(2y_{n})^{n}}f(y)$. The result follows. Q.E.D.

\ Now let us note that
\begin{equation}
D_{y}E(x,y)=(n-2)\frac{(\hat{x}-y)(\overline{x}-\overline{y})}{\|\hat{x}-y\|^{n}\|x-y\|^{n}}
\end{equation}
and

\begin{equation}D_{y}F(x,y)=(n-2)\frac{(x-y)(\overline{\hat{x}}-\overline{y})}{\|\hat{x}-y\|^{n}|x-y|^{n}}.
\end{equation}
Using these formulas we can deduce the following result.

\begin{theorem} Let $L\in C^{1}(\overline{K})$, then $I(y)$ is a hypermonogenic function on $R^{n}\backslash\overline{K}$ where

\[I(y)=y_{n}^{n-2}\left(\int_{K}E(x,y)L(x)dx^{n}-\int_{K}F(x,y)\hat{L}(x)dx^{n}\right).\]
\end{theorem}
{\bf{Proof:}}  First we have

\[\frac{n-2}{y_{n}}Q^{\prime}(I(y))=\frac{n-2}{y_{n}}\left[\frac{\hat{I}(y)-I(y)}{2}e_{n}\right]^{\prime}\]
\[=-\frac{n-2}{y_{n}}\left[\frac{\hat{I}(y)-I(y)}{2}\right]^{\prime}e_{n}=-\frac{n-2}{2y_{n}}e_{n}\left[I(y)-\hat{I}(y)\right]\]
By  (2) and (3) we have
\[M(I(y))=DI(y)+\frac{n-2}{y_{n}}Q^{\prime}(I(y))\]
\[=e_{n}(n-2)y_{n}^{n-3}\left[\int_{K} E(x,y) L(x)dx^{n}-\int_{K} F(x,y)\hat{L}(x)dx^{n}\right]\]
\[+y_{n}^{n-2}\left[\int_{K} D_{y}E(x,y) L(x)dx^{n}-\int_{K}D_{y }F(x,y)\hat{L(x)}dx^{n}\right]-\frac{(n-2)e_{n}}{2y_{n}}
\left[I(y)-\hat{I}(y)\right]\]
\[=\frac{(n-2)y_{n}^{n-3}}{2}\left[\int_{K}\frac{e_{n}(\overline{x}-\overline{y})}{\|x-y\|^{n}\|\hat{x}-y\|^{n-2}}L(x)dx^{n} \right.\]
\[+\int_{K}\frac{2y_{n}(\hat{x}-y)(\overline{x}-\overline{y})}{\|x-y\|^{n}\|\hat{x}-y\|^{n}}L(x)dx^{n}\]
\[\left. -\int_{K}e_{n}\frac{(\overline{x}-\overline{\hat{y}})}{\|x-y\|^{n-2}\|\hat{x}-y\|^{n}}L(x)dx^{n}\right]\]
\[+\frac{(n-2)y_{n}^{n-3}}{2}\left[\int_{K}\frac{-e_{n}(\overline{\hat{x}}-\overline{y})}{\|x-y\|^{n-2}\|\hat{x}-y\|^{n}}\hat{L}(x)dx^{n}\right. \]
\[-\int_{K}\frac{2y_{n}(x-y)(\overline{\hat{x}}-\overline{y})}{|x-y|^{n}\|\hat{x}-y\|^{n}}\hat{L}(x)dx^{n}\]
\[\left. +\int_{K} \frac{e_{n}(\overline{\hat{x}}-\overline{\hat{y}})}{\|x-y\|^{n}\|\hat{x}-y\|^{n-2}}\hat{L}(x)dx^{n}\right]\]
\[=\frac{(n-2)y_{n}^{n-3}}{2}\left[I_{1}+I_{2}\right].\]
Here
\[I_{1}=\int_{K}\frac{e_{n}(\overline{x}-\overline{y})\|\hat{x}-y\|^{2}+2y_{n}(\hat{x}-y)(\overline{x}-\overline{y})-e_{n}(\overline{x}-\overline{\hat{y}})\|x-y\|^2}{\|x-y\|^{n}\|\hat{x}-y\|^{n}}L(x)dx^{n}\]
and
\[e_{n}(\overline{x}-\overline{y})\|\hat{x}-y\|^{2}+2y_{n}(\hat{x}-y)(\overline{x}-\overline{y})-e_{n}(\overline{x}-\overline{\hat{y}})\|x-y\|^{2}\]
\[=\Big(e_{n}\|\hat{x}-y\|^{2}+2y_{n}(\hat{x}-y)-(\hat{x}-y)e_{n}(x-y)\Big)(\overline{x}-\overline{y})\]
\[=(\hat{x}-y)\Big((\overline{\hat{x}}-\overline{y})e_{n}-2y_{n}e_{n}e_{n}-(\overline{\hat{x}}-\overline{\hat{y}})e_{n}\Big)(\overline{x}-\overline{y})\]
\[=(\hat{x}-y)\Big((\overline{\hat{x}}-\overline{y})-2y_{n}e_{n}-(\overline{\hat{x}}-\overline{\hat{y}})\Big)e_{n}(\overline{x}-\overline{y})=0.\]

\ So $I_{1}=0$. Similarly

\[I_{2}=\int_{K}\frac{-\|x-y\|^{2}e_{n}(\overline{\hat{x}}-\overline{y})-2y_{n}(x-y)(\overline{\hat{x}}-\overline{y})+e_{n}(\hat{\overline{x}}-\hat{\overline{y}})\|\hat{x}-y\|^{2}}{\|x-y\|^{n}\|\hat{x}-y\|^{n}}\hat{L}(x)dx^{n}\]
and
\[-(x-y)(\overline{x}-\overline{y})e_{n}(\overline{\hat{x}}-\overline{y})-2y_n(x-y)(\overline{\hat{x}}-\overline{y})+(x-y)e_{n}\|\hat{x}-y\|^{2}\]
\[=(x-y)\Big(-(\overline{x}-\overline{y})(x-\hat{y})+2e_{n}y_{n}(x-\hat{y})+(\overline{x}-\overline{\hat{y}})(x-\hat{y})\Big)e_{n}\]
\[=(x-y)\Big(-\overline{x}+\overline{y}+2e_ny_n+\overline{x}-\overline{\hat{y}}\Big)(x-\hat{y})e_{n}=0.\]
So $I_{2}=0$ and $MI=0$. Q.E.D.

\begin{theorem}
 Let $F\in C^{1}(\overline{K})$, $y\in K$ then $M(I(y))=F(y)$ where $I(y)$ is as defined in Theorem 11.
\end{theorem}
{\bf{Proof:}}  First we show that $I(y)$ is a well defined  function on $K$. To do this we only need to show that the
integral  $$\int_{B(y,r)} E(x,y) F(x)dx^{n}-\int_{B(y,r)} F(x,y)\hat{L}(x)dx^{n}$$ is well defined for any ball
${B(y,r)}\subset K$. Since
$$\left|\left| \int_{B(y,r)} E(x,y) L(x)dx^{n}-\int_{B(y,r)} F(x,y)\hat{L}(x)dx^{n}\right| \right| \leq sup_{x\in
B(y,r)}\|L\|C(n)\int_{0}^{r} ds$$
the integral is clearly finite. Then we can calculate $MI(y)$. For any fixed point
$y\in K$ take any closed $n$-dimensional rectangle $R(y)\subset K$. Based on Theorem 13 we have
\[MI(y)=M\left[y_{n}^{n-2} \left(\int_{K}E(x,y) L(x)dx^{n}-\int_{K} F(x,y) \hat{L}(x)dx^{n}\right)\right]\]
\[=M\left[ y_{n}^{n-2} \left(\int_{K\backslash\overline{R}}E(x,y) L(x)dx^{n}-\int_{K\backslash\overline{R}} F(x,y)
\hat{L}(x)dx^{n} \right. \right. \]
\[\left. \left. +\int_{\overline{R}}E(x,y) L(x)dx^{n}-\int_{\overline{R}}F(x,y) \hat{L}(x)dx^{n}\right)\right]\]
\[=M\left[ y_{n}^{n-2} \left(\int_{\overline{R}} E(x,y) L(x)dx^{n}-\int_{\overline{R}}F(x,y) \hat{L}(x)dx^{n}\right)\right]\]

\ We consider

\[\lim\limits_{h_{j}\rightarrow 0}\frac{1}{h_j}\left[\int_{R(y)}\left(E(x,y-h_{j}e_{j})-E(x,y)\right)L(x)dx^{n}. \right. \]
\[\left. -\int_{R(y)}\left(F(x,y-h_{j}e_{j})-F(x,y)\right)\hat{L}(x)dx^{n}\right]\]
\[=\lim\limits_{h_{j}\rightarrow 0}\frac{1}{h_{j}}\left[\int_{R_{1}(y,h_{1})}E(x,y-h_{j}e_{j})L(x)dx^{n}\right. \]
\[+\int_{R(y,h_{j})\backslash R_{3}(y,h_{j})}F(x,y)\left(L(x-h_{j}e_{j})-L(x)\right)dx^{n}\]
\[-\int_{R_{3}(y,h_{j})}E(x,y)L(x)dx^{n} -\int_{R_{1}(y,h_{1})}F(x,y-h_{j}e_{j})\hat{L}(x)dx^{n}\]
\[+\int_{R(y,h_{j})\backslash R_{3}(y,h_{j})}F(x,y)\left(\hat{L}(x-h_{j}e_{j})-\hat{L}(x)\right)dx^{n}\]
\[\left. -\int_{R_{3}(y,h_{j})}F(x,y)\hat{L}(x)dx^{n}\right]\]
where $R_{1}(y,h_{j})$ is the closed rectangle obtained from $R(y)$ by truncating $R(y)$ in the $-e_{j}^{'}$th
direction a distance $h_{j}$ from the face whose normal vector is $-e_{j}$, $R_{2}(y,h_{j})=R(y)- R_{1}(y,h_{j})$ while
$R_{3}(y,h_{j})$ is the closed rectangle obtained from $R_{(y)}$ by truncating $R(y)$ in the$e_{j}^{'}$th direction a
distance $h_{j}$ from the face whose normal vector is $e_{j}$. The width of both $R_{1}(y,h_{j})$ and $R_{3}(y,h_{j})$
in the $e_{j}^{'}$th direction is $h_{j}$. Consequently the previous limits evaluates to
\[\frac{1}{h_{j}}\left[\int_{Q_{1}(y,j)}\left(E(x,y)L(x)-F(x,y)\right)\hat{L}(x)dx^{n}\right. \]
\[-\int_{Q_2(y,j)}\left(E(x,y)L(x)-F(x,y)\right)\hat{L}(x)dx^{n}\]
\[\left. +\int_{R(y)}E(x,y)\frac{\partial L(x)}{\partial x_{j}}dx^{n}-\int_{R(y)}F(x,y)
\frac{\partial \hat{L}(x)}{\partial x_{j}}dx^{n}\right].\] So $MI(y)=F(y)$ Q.E.D.

\section{M\"{o}bius Transformations and the Hyperbolic Dirac-Hodge Operator and Hyperbolic Laplacian}

\ We begin by establishing an invariance for the Cauchy Integral Formula under M\"{o}bius transformations. We begin by considering the case of Kelvin inversion $In(x)=-x^{-1}$ for $x\neq 0$. Suppose that $f(y)$ is left hypermonogenic on a domain $U$ in upper half space. For $K$ a closed bounded subregion of $U$ we have

\[f(y)=\frac{2^{n-2}y_{n}^{n-2}}{\omega_{n}}\int_{\partial K}\left(\frac{(x-y)}{\|x-y\|^{n}\|x-\hat{y}\|^{n-2}}n(x)f(x)\right. \]
\[\left. +\frac{(\hat{x}-y)}{\|\hat{x}-y\|^{n}\|x-y\|^{n-2}}\hat{n}(x)\hat{f}(x)\right)d\sigma(x)\]
for each $y$ in the interior of $K$. If now $y=-v^{-1}$ and $x=-u^{-1}$ then $y_{n}=\frac{v_{n}}{\|v\|^{2}}$ and the integral formula becomes
\[f(-v^{-1})=\frac{2^{n-2}v_{n}^{n-2}}{\| v\|^{2n-4}\omega_{n}}\int_{\partial K^{-1}}\left(v\|v\|^{2n-4}
\frac{(u-v)}{\|u-v\|^{n}\|\hat{u}-v\|^{n-2}}n(u)u^{-1}f(-u^{-1})\right. \]
\[\left. +v\|v\|^{2n-4}\frac{(\hat{u}-v)}{\|\hat{u}-v\|^{n}\|u-v\|^{n-2}}\hat{n}(u)\hat{u}^{-1}\hat{f}(-u^{-1})\right)d\sigma(u).\]
This expression simplifies to
\[v^{-1}f(-v^{-1})=\frac{2^{n-2}v_{n}^{n-2}}{\omega_{n}}\int_{\partial K^{-1}}\left(\frac{(u-v)}{\|u-v\|^{n}\|u-\hat{v}\|^{n-2}}
n(u)u^{-1}f(-u^{-1})\right. \]
\[\left. +\frac{(\hat{u}-v)}{\|\hat{u}-v\|^{n}\|u-v\|^{n-2}}\hat{n}(u)\hat{u}^{-1}\hat{f}(-u^{-1})\right)d\sigma(u).\]

\ Similar results may be obtained for the other generators of the M\"{o}bius group. It follows that we have:

\begin{theorem}
Suppose that $\psi(u)=(au+b)(cu+d)^{-1}$ is a M\"{o}bius transformation that leaves $R^{n,+}$ invariant. Suppose also that $f$ is left hypermonogenic on a domain $U\subset R^{n,+}$ and $K$ is a closed bounded subregion of $U$. Then
\[J(\psi,v)f\big(\psi(u)\big)=\frac{2^{n-2}v_{n}^{n-2}}{\omega_{n}}\int_{\partial\psi^{-1}(K)}
\left(\frac{(u-v)}{\|u-v\|^{n}\|\hat{u}-v\|^{n-2}}n(u)J(\psi,u)f(\psi(u))\right. \]
\[\left. +\frac{(\hat{u}-v)}{\|\hat{u}-v\|^{n}\|u-v\|^{n-2}}\hat{n}(u)\hat{J}(\psi,u)\hat{f}(\psi(u))\right)d\sigma(u)\]
where $J(\psi,u)=\frac{\widetilde{cu+d}}{\|cu+d\|^{2}}$.
\end{theorem}

\ Similarly one can take the function $I(y)$ set up in the statement of Theorem 14 and see that

\[I\big(\psi(v)\big)=\frac{v_{n}^{n-2}}{\omega_{n}\|v\|^{2n-4}}\int_{U^{-1}}\left(v^{-1}\|v\|^{2n-2}
\frac{(u-v)}{\|u-v\|^{n}\|\hat{u}-v\|^{n-2}}\frac{u}{\|u\|^{4}}L(\psi(u)) \right. \]
\[\left. -v^{-1}\|v\|^{2n-2}\frac{(\hat{u}-v)}{\|\hat{u}-v\|^{n}\|u-v\|^{n-2}}\frac{\hat{u}}{\|u\|^{4}}
\hat{L}(\psi(u))\right) du^{n}.\]
This simplifies to
\[v^{-1}I\big(\psi(v)\big)=\frac{v_{n}^{n-2}}{\omega_{n}}\int_{U^{-1}}\left(E(u,v)\frac{u}{\|u\|^{4}}L\big(\psi(u)\big)
-F(u,v)\frac{\hat{u}}{\|u\|^{4}}\hat{L}(u)\right)du^{n}.\]

\ If now we set $L(x)=M\phi(x)$ where $\phi$ has compact support in $U$ and apply the operator $M$ to the above
equation we obtain
\[M\left(v^{-1}f\Big(\phi\big(\psi(v)\big)\Big)\right)=\frac{v}{\|v\|^{4}}M\Big(\phi\big(\psi(v)\big)\Big).\]

\ Again similar results may be obtained for other generators of the conformal group. It follows that we have:
\begin{theorem}
Suppose that $\psi(u)=(au+b)(cu+d)^{-1}$ is a M\"{o}bius transformation that leaves $R^{n,+}$ invariant. Suppose also that $\phi$ is a $C^{1}$ function with support in $U$. Then
\[M\Big(J(\psi,v)\phi\big(\psi(v)\big)\Big)=J'(\psi,v)M\Big(\phi\big(\psi(v)\big)\Big)\]
where $J'(\psi,v)=\frac{\widetilde{cv+d}}{\|cv+d\|^{4}}$.
\end{theorem}

\ This theorem provides us with intertwining operators for the differential operator $M$ under actions of the
conformal group.\

 Using the previous theorem and a standard partition of unity argument we have:
\begin{proposition}
Suppose that $f:U\rightarrow Cl_{n}$ is a left hypermonogeinic function in the variable $x$ and $x=\psi(v)=(av+b)(cv+d)^{-1}$ is a M\"{o}bius transformation preserving upper half space then the function $J(\psi,v)f(\psi(v))$ is left hypermonogenic in the variable $v$.
\end{proposition}

\ This result was established in \cite{l} using different techniques.

\ Let us now consider the constant hypermonogenic function $f(x)=-e_{1}$. By the previous results then under inversion we obtain the left hypermonogenic function $-v^{-1}e_{1}=v'^{-1}$. The function $v'^{-1}$ is a direct analogue of the function $\frac{1}{z}$ from one complex variable. For each $k\in N$ the function $\frac{\partial^{k} v'^{-1}}{\partial v_{1}^{k}}=(-1)^{k}k!v'^{-k-1}$ is also left hypermonogenic on upper half space. Again by employing inversion it may now be observed that $v'^{k}$ is left hypermonogenic for each $k\in N$. These functions are direct analogues of the functions $z^{k}$ from one complex variable. That such functions are left hypermonogenic was first observed, using a different argument, in \cite{l}.

\ We will now proceed to find intertwining operators for the operators $\triangle_{R^{n,+}}$ and $\triangle'_{R^{n,+}}$ under M\"{o}bius transformations.

\ Using Proposition 2, Theorem 4 and (1) one may adapt the arguments used to establish  Theorem 16 to determine that for any M\"{o}bius transformation $\psi$ preserving upper half space and any $Cl_{n-1}$ valued $C^{2}$ function $\phi$ defined on a domain in upper half space
\[\triangle_{R^{n,+}}\Big(\phi\big(\psi(v)\big)\Big)=J_{1}(\psi,v)\triangle_{R^{n,+}}\Big(\phi\big(\psi(v)\big)\Big)\]
where $J_{1}(\psi,v)=\frac{1}{\|cv+d\|^{4}}$.

\ Similarly

\[\triangle'_{R^{n,+}}\Big(\phi\big(\psi(v)\big)\Big)=J_{1}(\psi,v)\triangle'_{R^{n,+}}\Big(\phi\big(\psi(v)\big)\Big).\]

\ It follows that if $\phi(x)$ is annihilated by $\triangle'_{R^{n,+}}$ then so is $\phi(\psi(v))$.


\begin{thebibliography}{99}
\bibitem{a} L. V. Ahlfors, {\it{M\"{o}bius Transformations in Several Dimensions}}, Ordway Lecture Notes, University of Minnesota, 1981.\bibitem{a1} L. V. Ahlfors, {\it{M\"{o}bius transformations in $R^{n}$ expressed through $2\times 2$ matrices of Clifford numbers}}, Complex Variables, 5, 1986, 215-224.
\bibitem{al} \"{O}. Akin and H. Leutwiler, {\it{On the invariance of the solutions of the Weinstein equation under M\"{o}bius transformations}}, K. Gowrisankran et al (eds), Classical and Modern Potential Theory and Applications, Kluwer, Dodrecht, 1994, 19-29.
\bibitem{ca} D. Calderbank, {\it{Dirac operators and Clifford analysis on manifolds}}, Max Plank Institute for Mathematics, Bonn, preprint number 96-131, 1996.
\bibitem{cw} C. Cao and P. Waterman, {\it{Conjugacy invariants of M\"{o}bius groups}}, Quasiconformal Mappings and Analysis (Ann Arbor, MI, 1995), Springer, New York, 1998, 109-139.
\bibitem{cc} P. Cerejeiras and J. Cnops, {\it{Hodge-Dirac operators for hyperbolic space}}, Complex Variables, 41, 2000, 267-278.
\bibitem{cn} J. Cnops, {\it{An Introduction to Dirac Operators on Manifolds}}, Progress in Mathematical Physics, Birkh\"{a}user, Boston, 2002.
\bibitem{e1} S.-L. Eriksson-Bique, {\it{M\"{o}bius transformations and $k$-hypermonogenic functions}} to appear.
\bibitem{e} S.-L. Eriksson, {\it{Integral formulas for hypermonogenic functions}}, to appear.
\bibitem{el2} S.-L. Eriksson-Bique {\it{$k$-hypermonogenic functions}}, Progress in Analysis, H. Begehr et al (editors), World Scientific, New Jersey, 2003, 337-348.
\bibitem{el1} S.-L. Eriksson-Bique and H. Leutwiler, {\it{Hypermonogenic functions}}, Clifford Algebras and their Applications in Mathematical Physics, Volume 2, ed J. Ryan and W. Spr\"{o}$\beta$ig, Birkh\"{a}user, Boston, 2000, 287-302.
\bibitem{el} S.-L. Eriksson and H. Leutwiler, {\it{Some integral formulas for hypermonogenic functions}}, to appear.
\bibitem{el3} S.-L. Eriksson and H. Leutwiler, {\it{Hypermonogenic functions and their Cauchy-type theorems}}, Trends in Mathematics: Advances in Analysis and Geometry, Birkh\"{a}user, Basel, 2003, 1-16.
\bibitem{glq} G. Gaudry, R. Long and T. Qian, {\it{A martingale proof of $L^{2}$-boundedness of Clifford valued singular integrals}}, Annali di Matematica, Pura Appl., 165, 1993, 369-394.
\bibitem{gs} K. Gowrisankran and D. Singman, {\it{Minimal fine limits for a class of potentials}}, Potential Anal., 13, 2000, 103-114.\bibitem{h} M. Habib, {\it{Invariance des fonctions $\alpha$-harmoniques par les transformations de M\"{o}bius}}, Exposition Math., 13, 1995, 469-480.
\bibitem{hu} A. Huber, {\it{On the uniqueness of generalized axially symmetric potentials}}, Ann. of Math., 60, 1954, 351-358.
\bibitem{i} V. Iftimie, {\it{Fonctions hypercomplexes}}, Bull. Math.. de la Soc. Sci. Math. de la R. S. de Roumanie, 9, 1965, 279-332.
\bibitem{kr} R. S. Krausshar and J. Ryan, {\it{Clifford and harmonic analysis on spheres and hyperbolas}} to appear in Revista Matematca Iberoamericana.
\bibitem{kr1} R. S. Krausshar, J. Ryan and Q. Yuying, {\it{Harmonic, monogenic and hypermonogenic functions on some conformally flat manifolds in $R^{n}$ arising from special arithmetic groups of the Vahlen group}}, to appear in Contemporary Mathematics.
\bibitem{kr2} R. S. Krausshar and J. Ryan, {\it{Some conformally flat spin manifolds, Dirac operators and automorphic forms}}, to appear.
\bibitem{lr} H. Liu, and J. Ryan, {\it{Clifford analysis techniques for spherical pde's}}, Journal of Fourier Analysis and its Applications, 8, 2002, 535-564.
\bibitem{l1} H. Leutwiler, {\it{Best constants in the Harnack inequality for the Weinstein equation}}, Aequationes Mathematicae, 34, 1987, 304-315.
\bibitem{l} H. Leutwiler, {\it{Modified Clifford analysis}}, Complex Variables, 17, 1992, 153-171.
\bibitem{Mc} A. McIntosh, {\it{Clifford algebras, Fourier theory, singular integrals, and harmonic functions on Lipschitz domains}}, Clifford Algebras in Analysis and Related Topics, J. Ryan (ed), CRC Press, Boca Raton, 1996, 33-87.
\bibitem{m} M. Mitrea, {\it{Singular Integrals, Hardy Spaces, and Clifford Wavelets}}, Lecture Notes in Mathematics, No 1575, Springer-Verlag, Heidelberg, 1994.
\bibitem{m1} M. Mitrea, {\it{Generalized Dirac operators on non-smooth manifolds and Maxwell's equations}}, Journal of Fourier Analysis and its Applications, 7, 2001, 207-256.
\bibitem{p} I. Porteous, {\it{Clifford Algebras and the Classical Groups}}, Cambridge University Press, Cambridge, 1995.
\bibitem{r} J. Ryan, {\it{Dirac operators on spheres and hyperbolae}}, Bolletin de la Sociedad Matematica a Mexicana, 3, 1996, 255-270.
\bibitem{v} K. Th. Vahlen, {\it{\"{U}ber Bewegungen und Complexe Zahlen}}, Math. Ann., 55, 1902, 585-593.
\bibitem{vl} P. Van Lancker, {\it{Clifford analysis on the sphere}}, Clifford Algebras and their Applications in Mathematical Physics, V. Dietrich et al (editors), Kluwer, Dordrecht, 1998, 201-215.
\bibitem{w} A. Weinstein, {\it{Generalized axially symmetric potential theory}}, Bull. Amer. Math. Soc., 59, 1953, 20-38.
\end{thebibliography}
\end{document}